\title{Reduced basis approaches for parametrized bifurcation problems held by non-linear {V}on {K{\'a}rm{\'a}n} equations}
\author{
        Federico Pichi\\
        mathLab, Mathematics Area, SISSA,\\ International School for Advanced Studies \\           	    via Bonomea 265, 34136 Trieste, Italy\\
        \and
        Gianluigi Rozza \\
        mathLab, Mathematics Area, SISSA, \\International School for Advanced Studies \\           	    via Bonomea 265, 34136 Trieste, Italy\\
}
\date{\today}

\documentclass[12pt]{article}
\usepackage{mathptmx}       % selects Times Roman as basic font
\usepackage{helvet}         % selects Helvetica as sans-serif font
\usepackage{courier}        % selects Courier as typewriter font
\usepackage{type1cm}        % activate if the above 3 fonts are
\usepackage{makeidx}         % allows index generation
\usepackage{graphicx}        % standard LaTeX graphics tool
\usepackage{multicol}        % used for the two-column index
\usepackage[bottom]{footmisc}% places footnotes at page bottom

\usepackage{amssymb}
\usepackage{amsmath}
\usepackage{graphicx}
\usepackage{epsfig}
\usepackage{subfigure}
\usepackage{multirow}
\usepackage{color}
\usepackage{newtxtext,newtxmath}
\usepackage{graphicx}
\usepackage{epsfig}
\usepackage{subfigure}
\usepackage{multirow}
\usepackage{color}
\usepackage{breqn}
\usepackage{listings}
\usepackage{xcolor}
\usepackage{booktabs}
\usepackage{siunitx}
\usepackage{xcolor}
\usepackage{colortbl}
\usepackage{url}
\usepackage{algorithm}
\usepackage{algpseudocode}
\usepackage{float}

\definecolor{light-gray}{gray}{0.95}

\graphicspath{{Figures/}}

\newcommand{\p}{\partial}
\graphicspath{{Figures/}}

\newcommand{\R}{\mathbb{R}}
\newcommand{\lam}{\lambda}
\newcommand{\Om}{\Omega}
\newcommand{\D}{\mathcal{D}}
\newcommand{\vK}{\text{{V}on {K{\'a}rm{\'a}n}}}
\newcommand{\ackname}{Acknowledgements}
\newtheorem{definition}{Definition}
\newtheorem{theorem}{Theorem}

\begin{document}

\title{Reduced basis approaches for parametrized bifurcation problems held by non-linear {V}on {K{\'a}rm{\'a}n} equations}

\maketitle

\begin{abstract}
This work focuses on the computationally efficient detection of the buckling phenomena and bifurcation analysis of the parametric {V}on {K{\'a}rm{\'a}n} plate equations based on reduced order methods and spectral analysis. The computational complexity - due to the fourth order derivative terms, the non-linearity and the parameter dependence - provides an interesting benchmark to test the importance of the reduction strategies, during the construction of the bifurcation diagram by varying the parameter(s). To this end, together the state equations, we carry out also an analysis of the linearized eigenvalue problem, that allows us to better understand the physical behaviour near the bifurcation points, where we lose the uniqueness of solution. We test this automatic methodology also in the two parameter case, understanding the evolution of the first buckling mode.
\end{abstract}

\section{Introduction and motivation}
\label{sec:1}
In this work we are interested in the numerical approximation of parameter dependent non-linear structural problems governed by the $\vK$ plate equations \cite{vonka}. The main issue, and so the most interesting feature, of this kind of model, is the bifurcation phenomenon that corresponds to the buckling behaviour of an elastic thin plate.

 From the mathematical point of view, we dealt with a problem in which multiple solutions for the same value of the parameter can arise. This led us to the non-uniqueness of the solution that we investigated. Numerically speaking the model presents three main problems that we had to face with: (i) non-linearity (ii) parameter dependency (iii) high order derivative terms. 
 
 Thus we relied on the Reduced Basis (RB) method aiming at reducing the computational time in order to have a better understanding of the physical phenomena that we were modelling. In practice we wanted to find an efficient and rapid way to draw the bifurcation diagram and detect the buckling points, i.e. critical values of the parameter $\lam$, which in our case controls the compression along the edges of the plate.

The RB method is one of the main Reduced Order Modelling (ROM) techniques and we applied it in conjunction with the Galerkin-Finite Element (FE) method \cite{patera07:book,hesthaven2015certified,quarteroni2015reduced,rozza08:ARCME}. As the latter, the RB approach is a Galerkin projection over a finite dimensional subspace of the weak formulation of the model. Thus, we solve the FE problem, called high order formulation because of the huge number of degrees of freedom involved, then we construct a subspace as the span of some basis functions computed before, and finally we project again on this new smaller subspace. A recent review chapter focused on parametric elasticity problems solved by RB method is \cite{Huynh2018}.

This technique allows us to efficiently study the entire behaviour of the plate under compression, but still involves a huge amount of computations in the offline high order phase. Thus, being inspired by the recent works in branching detection \cite{herrero2013rb,pitton2017}, also in a more industrial framework \cite{ZAMM:ZAMM201500217}, we supplemented the model with the eigenvalue analysis of the linearized equations.

Indeed, in the past, many authors established the connection between the buckling points and the eigenvalues behaviour, but the computational complexity of the problem was not affordable for that time \cite{bauerreiss,berger}. Hence we propose a new reduced order approach, in order to avoid the huge computational cost.  Previous works on model order detection for non-linear elasticity could be found in \cite{veroy03:_phd_thesis,zanon:_phd_thesis}, as well as in preliminary works by Noor and Peters \cite{noor80:_reduc,noor82,noor83:_reduc}.

The structure of the work is the following. In Section $2$ we provide a brief description of the equations of plate $(2.1)$, that model the physical phenomena in connection with the mathematical formulation and boundary conditions, then we provide the weak formulation $(2.2)$ of the problem, necessary as the first step towards the numerical approximation. To end, we recall few definitions and theorems in bifurcation and non-linear analysis $(2.3)$ that justify the eigenproblem coupling.

In Section $3$  we deal with the numerical approximation of the problem, providing the pseudo-code used to construct the bifurcation diagram, based on Galerkin finite element $(3.1)$ and Newton method. The reduction strategy, or reduced basis method is reported in $(3.2)$ with its main features and its matrix formulation, compared with the finite element one. Finally, we show some preliminary results on the spectral analysis $(3.3)$, the eigenvalues approximation in two different settings, and a test to verify the order of convergence which is fundamental in view of the connection with the buckling points.

Section $4$ is dedicated to results and tests. Here we ensure the reliability of our high fidelity solver and most importantly we ran all the procedure to study the square and rectangular plate case, finding up to eight solutions for the same parameter value, multiple buckling points that validate the theoretical results, and provide a good accuracy for the reduced approximation, while saving significant computational time. 
Finally we show some preliminary results in the two-parameters test case, where the shape of the compression load is also parametrized, showing a 3-D bifurcation plot for the evolution of the first buckling mode, which is very complex and computationally very expansive. To the best of our knowledge the proposed approach combining reduced order methods and parametric bifurcation analysis is original, especially considering more then one parameter for non-linear $\vK$ equations. Some conclusions follow.

\section{Parametrized formulation of {V}on {K{\'a}rm{\'a}n} equations}
Starting from the very well known theory of continuum mechanics, {V}on {K{\'a}rm{\'a}n} in 1910 proposed a mathematical model in order to describe all the possible configurations that a plate under compression can take \cite{vonka}.
\emph{Buckling phenomenon} is the mathematical way of explaining a well known physical event that very frequently happens in many contexts. As an example, an appropriate one since it is exactly what we are trying to model, if we pick a thin rectangular plate at rest, we can use our hands to compress it until we reach a critical point, i.e. when the plate takes a deformed configuration, or it buckles.  

\subsection{Equations of the plate}

Let us consider an elastic and rectangular plate $\Omega = [0,L] \times [0,1]$ in its undeformed state, subject to a $\lambda$-parametrized external load acting on its edges, then the displacement from its flat state and the \emph{Airy stress potential}, respectively $u$ and $\phi$, satisfy the {V}on {K{\'a}rm{\'a}n} equations
\begin{equation}
\label{eq:karm}
\begin{cases}
\Delta^2u =  \left[\lam h + \phi, u\right] + f \ , \quad &\text{in}\ \Omega \\
\Delta^2\phi = -\left[u, u\right] \ , \quad &\text{in}\ \Omega 
\end{cases}
\end{equation}
where $h$ and $f$ are some given functions, that we can set to specify the external forces acting on our plate, while 
%$$\Delta^2 := \Delta\Delta = \left(\frac{\p\,^2}{\p\,x^2} + \frac{\p\,^2}{\p\,y^2}\right)^2 \ ,$$
$\Delta^2$ is the biharmonic operator in Cartesian coordinates and 
$$[u,\phi] := \frac{\p\,^2u}{\p\,x^2}\frac{\p\,^2\phi}{\p\,y^2} -2\frac{\p\,^2u}{\p x\p y}\frac{\p\,^2\phi}{\p x\p y} + \frac{\p\,^2u}{\p\,y^2}\frac{\p\,^2\phi}{\p\,x^2} \ ,$$
is the \emph{brackets of Monge-Amp{\'e}re}.
 Thus we aim to find the displacement and the coupled \emph{Airy stress potential}
 % (physically linked to the second derivatives of the Piola-Kirchhoff stress tensor) 
 that solve the system \eqref{eq:karm} which is of fourth order, due to the presence of the biharmonic operator, non-linear due to the product of second derivatives in the bracket, and parametric due to the buckling coefficient $\lam$ varying in a proper range of real numbers. Moreover we are presenting, for the sake of simplicity, a non-dimensional model where all the physical quantities, except for the compression parameter $\lam$,  are set to unity.
 
In order to have a well posed system of partial differential equations we provide boundary conditions for both the unknowns that match the different physical setting. Among all the possible choices, we just focused on the so called 
%We present only the most used ones, but they can be imposed in many different ways in order to be coherent with the experiments of interest \cite{ciarlet1997mathematical}. The first option is to impose \textit{totally clamped} boundary conditions of the form
%\begin{equation}
%\label{eq:bctc}
%\begin{cases}
%u = \p_n u =  0 , \quad &\text{in}\ \p\Omega \\
%\phi = \p_n \phi =  0 , \quad &\text{in}\ \p\Omega 
%\end{cases}
%\end{equation}
%whose meaning is that the plate is completely blocked on its sides. Note that we denoted with $\p_n$ the directional derivative along the normal $n$ to the boundary $\p \Om$. 
%A second possible option is the so called 
\textit{simply supported} boundary conditions
\begin{equation}
\label{eq:bcss}
\begin{cases}
u = \Delta u =  0 , \quad &\text{in}\ \p\Omega \\
\phi = \Delta \phi =  0 , \quad &\text{in}\ \p\Omega 
\end{cases}
\end{equation}
which are physically complex to reproduce, but also the most used ones for the simulations because of their versatility and importance also in the weak formulation. So from now on, unless specified otherwise, we will consider the system \eqref{eq:karm} with simply supported boundary conditions \eqref{eq:bcss}. 
We remark that despite the simple boundary conditions chosen, the goal of this work is understanding the bifurcation behaviour for a complex system, regardless the numerical constraints that a conforming method for more involved boundary condition could impose.

%Before moving on to the weak formulation let us consider another useful way of expressing the {V}on {K{\'a}rm{\'a}n} equations. In fact it seems to be natural, once chosen the boundary conditions \eqref{eq:bcss},
Thanks to the BCs chosen we can split the system of two fourth order non-linear elliptic equations, into a system of four second order non-linear elliptic equations.
In order to carry out this trick, given by Ciarlet-Raviart \cite{ciarlet1974mixed}, we introduce two new unknowns, namely $U = \Delta u, \ \Phi = \Delta \phi$, so that we can rewrite the system \eqref{eq:karm} with homogeneous Dirichlet boundary conditions as 
\begin{equation}
\label{eq:karmsplit}
\begin{cases}
\Delta U =  \left[\lam h + \phi, u\right] + f \ , \quad &\text{in}\ \Omega \\
\Delta u = U \ , \quad &\text{in}\ \Omega \\
\Delta \Phi = -\left[u, u\right] \ , \quad &\text{in}\ \Omega \\
\Delta \phi = \Phi  \ , \quad &\text{in}\ \Omega 
\end{cases} \qquad \text{with} \qquad \begin{cases} u =  0 , \quad &\text{in}\ \p\Omega \\
U =  0 , \quad &\text{in}\ \p\Omega \\
\phi = 0 , \quad &\text{in}\ \p\Omega \\
\Phi =  0 . \quad &\text{in}\ \p\Omega
\end{cases}
\end{equation}
%to which we assign the homogeneous Dirichlet boundary conditions, derived from the simply supported ones
%\begin{equation}
%\label{eq:bcsssplit}
%\begin{cases}
%u =  0 , \quad &\text{in}\ \p\Omega \\
%U =  0 , \quad &\text{in}\ \p\Omega \\
%\phi = 0 , \quad &\text{in}\ \p\Omega \\
%\Phi =  0 . \quad &\text{in}\ \p\Omega
%\end{cases}
%\end{equation}

We know that \eqref{eq:karm} and \eqref{eq:karmsplit} are equivalent \cite{zhang2008invalidity} when the boundary is regular and the solution is smooth enough, so from now on we just consider the latter.

\begin{figure} [htbp]   
\centering    
\input{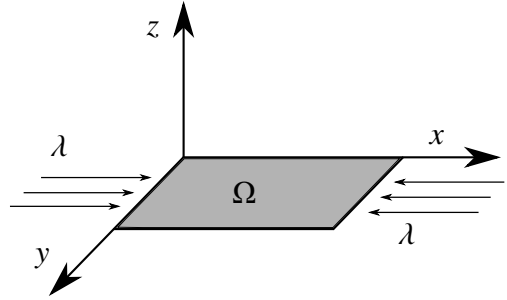}   
\caption{A rectangular bi-dimensional elastic plate compressed on its edges}  
\label{fig:plate} 
\end{figure}

Finally, since we are interested in the behaviour of the plate under compression, we can set the external body force $f = 0$ 
%and study the homogeneous system. Moreover through the function $h$ we can 
and model different kind of stresses at the boundaries through the function $h$. Indeed if we choose $h = -\frac{1}{2}y^2$, we obtain that $\left[\lam h , u\right] = -\lam u_{xx}$
%system reads as
%\begin{equation}
%\label{eq:karmsplittheta}
%\begin{cases}
%\Delta U + \lam u_{xx}=  \left[\phi, u\right] \ , \quad &\text{in}\ \Omega \\
%\Delta u = U \ , \quad &\text{in}\ \Omega \\
%\Delta \Phi = -\left[u, u\right] \ , \quad &\text{in}\ \Omega \\
%\Delta \phi = \Phi  \ , \quad &\text{in}\ \Omega 
%\end{cases}
%\end{equation}
where we are assuming that the compression is acting on the edges parallel to the $y$ direction (see Figure \ref{fig:plate}).
Note also that if instead we choose $h = -\frac{1}{2}(x^2 + y^2)$ we would have the stress component given by $\left[\lam h , u\right] = -\lam \Delta u$, in which case the compression we are considering is on the whole boundary.

\subsection{Weak formulation}
Starting from the considerations of the previous section we are now able to set our problem in a more abstract mathematical framework, which will be considered by us for the following numerical investigation.
So let us consider $\lam \in \D \subset \mathbb{R}$, where $\D$ is the parameter space, $\Om \subset \mathbb{R}^2$ is the bi-dimensional domain, that we identify with the plate, whereas $V = V(\Om) =  \left(H^1_0(\Om)\right)^4$ is the Hilbert space in which we will seek the solution and $V'$ its dual space. 

Furthermore we represent the non-linear PDE with the parametrized mapping $G: V \times \D \rightarrow V'$, so that the {V}on {K{\'a}rm{\'a}n} system \eqref{eq:karmsplit}, in abstract form reads: given $\lam \in \D$, find $X(\lam) \doteq (u(\lam), U(\lam), \phi(\lam), \Phi(\lam)) \in V$ such that 
\begin{equation}
\label{eq:abstractform}
G(X(\lam); \lam) = 0\ , \quad \text{in} \ V'. 
\end{equation}

We present the weak formulation where all the boundary terms vanish due to the simply supported boundary conditions \eqref{eq:bcss}. Then, we seek $X(\lam) \in V$  such that
\begin{equation}
\label{eq:weakform}
\begin{cases}
(\nabla u, \nabla w)_{L^2(\Om)} + (U,w)_{L^2(\Om)} = 0 \ , \quad & \forall\,  w \in H_0^1(\Om)  \\
(\nabla U, \nabla v)_{L^2(\Om)} + (\left[\lam h + \phi, u\right], v)_{L^2(\Om)} = 0 \ , \quad & \forall\,  v \in H_0^1(\Om)  \\
(\nabla \phi, \nabla \theta)_{L^2(\Om)} + (\Phi,\theta)_{L^2(\Om)} = 0 \ , \quad & \forall\,  \theta \in H_0^1(\Om)  \\
(\nabla \Phi, \nabla \psi)_{L^2(\Om)} - (\left[u, u\right], \psi)_{L^2(\Om)} = 0 \ , \quad & \forall\,  \psi \in H_0^1(\Om)  
\end{cases}
\end{equation}
in which we embed the simply supported boundary conditions in the choice of the space $H^1_0(\Om)$, where the test functions $Y \doteq (w, v, \theta, \psi)$ reside. Moreover we denote with $(\cdot , \cdot)_{L^2(\Om)}$ the usual inner product in the Hilbert space $L^2(\Om)$.

Then coming back to the abstract form of our problem \eqref{eq:abstractform} the weak formulation reads: given $\lam \in \D$, find $X(\lam) \in V$ such that
\begin{equation}
\label{eq:abstractweakform}
g(X(\lam),Y ; \lam) = 0\ , \quad \forall\ Y \in V, 
\end{equation}
where the parametrized variational form $g(\cdot , \cdot ; \lam) : V \times V \rightarrow \mathbb{R}$ is defined as
%\begin{equation}
%\label{eq:weakdef}
$g(Z , Y ; \lam) = \langle G(Z ; \lam), Y\rangle\ , \ \forall\ Z, Y \in V,$
%\end{equation}
where we denoted the duality pairing between $V'$ and $V$ with $\langle \cdot, \cdot \rangle$.
In this case the parametrized variational form $g(\cdot , \cdot ; \lam)$ is defined as follows
%\begin{equation}
%\label{eq:vkweak}
%\begin{split}
%g(X(\lam),Y ; \lam) =& \ a(u(\lam),w) + b(U(\lam),w) + a(U(\lam),v) + \lam c(h, u(\lam), v) + c(\phi(\lam), u(\lam), v)\,  + \\
%& \ a(\phi(\lam), \theta) + b(\Phi(\lam), \theta) + a(\phi(\lam), \psi) - c(u(\lam), u(\lam), \psi)\ , \quad \forall\ Y \in V, \ \forall \lam \in \D
%\end{split}
%\end{equation}
\begin{dmath}
g(X(\lam),Y ; \lam) = \ a(u(\lam),w) + b(U(\lam),w) + a(U(\lam),v) + \lam c(h, u(\lam), v) + c(\phi(\lam), u(\lam), v) + a(\phi(\lam), \theta) + b(\Phi(\lam), \theta) + a(\Phi(\lam), \psi) - c(u(\lam), u(\lam), \psi)\ , {\quad \forall\ Y \in V, \ \forall \lam \in \D , }
\label{eq:vkweak}
\end{dmath}
where the following bilinear and trilinear forms have been introduced
\begin{equation*}
a(x, y) = \int_\Om{\nabla x \cdot \nabla y \ d\Om} \ , \qquad
b(x, y) = \int_\Om{xy \ d\Om} \ , \qquad
c(x, y, z) = \int_\Om{\left[x, y\right]z \ d\Om} \ .
\end{equation*}
The numerical treatment of the variational form including the bracket of Monge-Amp{\'e}re obviously needs a non-linear method, to this end we compute the partial Fr{\'e}chet derivative of $g(Z, \cdot ; \lam)$ with respect to $X$ at $Z \in V$.

 Thus we assume the mapping $G$ to be continuously differentiable and denote by $D_XG(Z; \lam) : V \rightarrow V'$ its partial Fr{\'e}chet derivative at $(Z, \lam) \in V \times \D$. In this way we can express the partial Fr{\'e}chet derivative of $g(Z, \cdot ; \lam)$ at $Z \in V$ as 
\begin{equation}
\label{eq:frechetder}
dg[Z](W, Y; \lam) = \langle D_XG(Z; \lam)W, Y \rangle , \quad \forall \ W,Y \in V.
\end{equation}

These computations, for the non-linear system we are considering, show the explicit expression of the derivative of $g$ be of the form
\begin{dmath}
\label{eq:frechetdervk}
dg[Z](X(\lam), Y; \lam) =  \ a(u(\lam),w) + b(U(\lam),w) + a(U(\lam),v) + \lam c(h, u(\lam), v) + c(\phi(\lam), Z_1, v) + c(Z_3, u(\lam), v)  + a(\phi(\lam), \theta) + b(\Phi(\lam), \theta) + a(\Phi(\lam), \psi) - c(u(\lam), Z_1, \psi) - c(Z_1, u(\lam), \psi)\ , {\quad \forall\ Z,Y \in V, \ \forall \lam \in \D , }
\end{dmath}
where we denoted with  $Z = (Z_1, Z_2, Z_3, Z_4)$ the components of the point in which we are computing the derivative.

\subsection{Bifurcation and non-linear analysis}
The focus of our work is the efficient detection of the possible multiple solutions of the equations \eqref{eq:karmsplit}. 
%In order to better explain this lack of uniqueness, let us go back to the physics behind the $\vK$ model. 
%
%We note that, due to the symmetry of the problem, when the plate starts buckling we can expect at least two different configurations. In fact if a given displacement is a solution of the system for a given value of the compression parameter $\lam$, then of course the displacement that we can obtain reflecting this configuration with respect to the plane where the plate lies, is again a solution of the system. Moreover, for each value of $\lam \in \R$, the equations \eqref{eq:karm} with zero external force (namely $f = 0$) always admit the null solution $(u, \phi) = (0, 0)$, which corresponds to the original undeformed configuration. 
%
%From these considerations, we understand that if we select a value for the compression, then the problem can be ill-posed: which solution should we expect from our numerical solver? How can we know that the solution we found is the unique one?
Following the work done in \cite{ambrosetti1995primer,caloz1997numerical,Brezzi1980} we recall the mathematical definitions of bifurcation theory, which, as we will see, will serve us in developing a tool for the efficient detection of the buckling points.

%Let us consider again the functional equation \eqref{eq:abstractform}. 
%As we noted also from an empirical point of view, 
Since the undeformed configuration is a \textit{trivial solution} for every $\lam \in \R$, i.e. $G(0, \lam) = 0$, we can denote with
%, which will be referred to us the \textit{trivial solution}. Moreover if we denote with 
\begin{equation}
\label{eq:solset}
\mathcal{S} = \{(X,\lam) \in V \times \R : X \neq 0 ,\  G(X, \lam) = 0\} \ ,
\end{equation}
the set of \textit{non-trivial} solutions of \eqref{eq:abstractform}, then we can finally define the bifurcation points.

\begin{definition}
We say that $\lam^* \in \R$ is a bifurcation point for $G: V \times \D \rightarrow V'$, from the trivial solution, if there is a sequence $(X_n,\lam_n) \in V \times \R$ with $X_n \neq 0$ and $G(X_n, \lam_n) = 0$ such that $$(X_n, \lam_n) \rightarrow (0, \lam^*)\ .$$
\end{definition}
%Equivalently, we are requiring the point $(0, \lam^*)$ to belong to the closure in $V \times \R$ of $\mathcal{S}$, i.e. in any neighbourhood of $(0, \lam^*)$ there is a point $(X, \lam) \in \mathcal{S}$. 
%Now the question is: how can we investigate in an efficient way the manifold $\mathcal{S}$ to understand where are these bifurcation points located? 

In order to understand where are located these bifurcation points, we could numerically investigate the equations for each value of the parameter $\lam$ observing when the buckling phenomena occur. Of course this way is computationally very expensive and we need a different tool for the detection. 

To this end we have analyzed the path pursued in \cite{berger,ambrosetti1995primer,bauerreiss} where the link between the bifurcation points and the behaviour of the eigenvalues of the linearized problem is highlighted. If we linearize the equations \eqref{eq:karmsplit} around the trivial solution the system we obtain is simply given by  
\begin{equation}
\label{eq:linkarmsplit}
\begin{cases}
\Delta U = \lam \left[h, u\right] \ , \quad &\text{in}\ \Omega \\
\Delta u = U \ , \quad &\text{in}\ \Omega 
\end{cases}, \qquad \text{with} \qquad
%\end{equation}
%with the homogeneous boundary conditions 
%\begin{equation}
%\label{eq:linbcsssplit}
\begin{cases}
u =  0 , \quad &\text{in}\ \p\Omega \\
U =  0 , \quad &\text{in}\ \p\Omega \, .
\end{cases} 
\end{equation}

This connection is not surprising, in fact from ODE's theory we know that the stability of the solutions is linked to eigenvalues that change sign, i.e. cross the imaginary axis varying $\lam$.
%we know that to study the stability of an ODE we search for the sign properties of the eigenvalues. Thus what we should expect in this framework is the parametrized eigenvalue of the linearized problem, connected to the bifurcation, that  crosses the imaginary axis varying $\lam$.
Now we briefly recall the main theorems \cite{ambrosetti1995primer}, based on an application of the Implicit Function Theorem, that validate the numerical investigation we carried out. 
\begin{theorem}
\label{th:ambro}
A necessary condition for $\lam^*$ to be a bifurcation point for $G: V \times \D \rightarrow V'$ is that the partial derivative $D_XG(0; \lam^*)$ is not invertible.
\end{theorem}
%This is a general result valid for an abstract equation such as the one in \eqref{eq:abstractform}, furthermore if we focus on $\vK$ equations we have the following theorem \cite{berger,bauerreiss}.
\begin{theorem}
\label{berg1}
Bifurcation points of $\vK$ equations \eqref{eq:karmsplit} with respect to the trivial solution, i.e. $X = 0$, can occur only at the eigenvalues of the linearized problem \eqref{eq:linkarmsplit}. 
\end{theorem}
The former is a general result, while for $\vK$ equations we also know \cite{berger,bauerreiss} that every bifurcation point is an eigenvalue of the linearized problem. %, but we also know that in this case the assertion can be inverted. 
Moreover, if we assume that all the eigenvalues are real, positive and ordered in such a way that $0 \leq \lam_1 \leq \lam_2 \leq \lam_3 \leq \ldots$, from \cite{berger} we can reverse the statement.%have the following validation to the numerical results we will show later.
\begin{theorem}
\label{berg3}
From each eigenvalue of the system \eqref{eq:linkarmsplit} at least one branch of non-trivial solution of \eqref{eq:karm} bifurcates. In particular from a simple eigenvalue one branch bifurcates and from a multiple eigenvalue at least two branches bifurcate. Furthermore if $\lam_1$ is the smallest eigenvalue then for every value $\lam \leq \lam_1$ the unforced system has no non-trivial solutions.
\end{theorem}

%Thus we have fully presented the relation between the bifurcation points and the spectral analysis of the linearized problem. 
%
%Now we have all the definition needed to set the $\vK$ model in finite dimension, that is finding the numerical solution that efficiently approximates the real one, investigating the bifurcation diagram and of course exploring the eigenvalues properties of the system. 
%
\section{Numerical approximation of the problem}
As we said, in this application
%noted in the previous section, the equations we want to solve present multiple ``critical points'', in the sense that 
we have to face with different kinds of difficulties, such as: non-linearity, fourth order derivatives, and parameter dependence.

This led us to use a double approach, the \emph{full order} and the \emph{reduced order} one that we will present soon in detail. 
Let us present now, through a pseudo-code, how we addressed all these issues in order to obtain the bifurcation diagram.
%Now we want to give an idea, through a pseudo-code, how  we addressed all these issues in order to obtain the bifurcation diagram. l
%\vspace{0.3cm}
\begin{algorithm}
\caption{A pseudo-code setting}\label{alg:01}
\begin{algorithmic}[1]
\While{$\lambda < \lambda_{f}$}\Comment{External loop on compression parameter}
	\If{$||u_h||_\infty < \delta$} \Comment{Pre-buckling initial guess}
		\State{$X_h(\lam) = X_h^{(0)}$}
	\EndIf
	\While{$||\delta X_h||_{V} > \epsilon$}\Comment{Newton method}
		\State{$D_XG(X_h^{(i)}(\lam); \lam)\,\delta X_h = G(X_h^{(i)}(\lam); \lam)$}\Comment{Galerkin FE method}
		\State{$X_h^{(i+1)}(\lam) = X_h^{(i)}(\lam) - \delta X_h$}
	\EndWhile
	\State{$\lambda = \lambda + d\lambda$}
\EndWhile
\end{algorithmic}
\end{algorithm}
%\begin{lstlisting}
%while ($\lambda < \lambda_{f}$)
%
% if($||u||_\infty < \delta$)
%    $X(\lam) = X^{(0)}$;
% end	
%
% while($||\delta X||_{V} > \epsilon$)
%    $D_XG(X^{(i)}(\lam); \lam)\delta X = G(X^{(i)}(\lam); \lam)$;
%    $X^{(i+1)}(\lam) = X^{(i)}(\lam) - \delta X$;
% end
%
% $\lambda = \lambda + d\lambda$;
%
%end
%\end{lstlisting}

The Algorithm \ref{alg:01} is the implementation result of our approach to the difficulties previously mentioned. We denote by $X_h = (u_h, U_h, \phi_h, \Phi_h)$ the approximated solution and by $ \D = [\lam_i, \lam_f ]$ the parameter domain
%mainly based on three fundamental parts which correspond to the three methods we implemented. 
and we start with a slightly modified \textit{continuation method} \cite{allogwer}, which in its basic formulation consists in a while loop over $\lambda \in \D$. So, at each new cycle, i.e. for every $\lam$, we are changing the initial guess for the non-linear problem from a properly chosen guess (usually the solution of the linearized problem) to the solution of the last cycle, as soon as the solution of the latter turns out to be non-trivial. Thus we are able to detect the buckling by looking at the behaviour of the maximum norm of the displacement.

%This happens when we reach the critical $\lambda$ value, or bifurcation point. Thus this is the key point in order to study how the solution varies with respect to the parameter, and so this is the tool we used to draw the bifurcation diagram.
The second goal we had was to overcome the issue of the non-linearity.
% that was related to find the solution of a problem of the type \eqref{eq:abstractform}. 
We choose the very well know \emph{Newton-Kantorovich method} \cite{ciarlet2013linear}, where we use the norm of the increment $\delta X_h$ in the Hilbert space $V$ as stopping criterion.

Finally we note that in the line 6 we have to find the solution of a new, linear, weak formulation. To this end we applied the \textit{Galerkin-Finite Element  method} which appears to be a good candidate, also in view of the numerical extension towards the model order reduction. 
Since the fundamental importance of the last two methods, we briefly recall in the next section how we applied them to the $\vK$ equations.
% This will serve us also to introduce the reduction strategies.

\subsection{Galerkin finite element method}
As we already said, we decided to use the Galerkin finite element method to discretize the problem. This is a projection-like technique, where using the versatility of the weak formulation, we can set the model in a finite dimensional space \cite{ciarlet2002finite,quarteroni2017numerical}.

So we again consider the weak problem \eqref{eq:abstractweakform} and denote with $V_h$ a family of spaces, dependent from the $h$ parameter, such that $V_h \subset V$, $dim(V_h) = N_h < \infty$ for all $h > 0$. 
%We are modelling a rectangular plate, so the discretized domain $\Om_h = int \left(\bigcup_{K \in \mathcal{T}_h} K \right)$, i.e. the interior of the union of the elements of the triangulation $\mathcal{T}_h$, coincides exactly with the real domain $\Om$. 
Then, given $\lam \in \D$, we seek $X_h(\lam) \in V_h$ that satisfies 
\begin{equation}
\label{eq:weakgal}
g(X_h(\lam),Y_h ; \lam) = 0\ , \quad \forall\ Y_h \in V_h. 
\end{equation}
Obviously, because of the non-lineartity, we can not directly apply the finite element method to the equation \eqref{eq:weakgal}. So the Newton method, which in this case is also known as Newton-Kantorovich  \cite{ciarlet2013linear,quarteroni2015reduced}, reads as follows: once assigned an initial guess $X^0_h(\lam) \in V_h$, for every $k = 0, 1, \dots$ we seek the variation $\delta X_h \in V_h$ such that
\begin{equation}
\label{eq:weaknewton}
dg[X_h^k(\lam)](\delta X_h, Y_h; \lam) =  g(X_h^k(\lam),Y_h ; \lam) \ , \quad \forall\ Y_h \in V_h, 
\end{equation}
and then we update the solution for the successive step as $X_h^{k+1}(\lam) = X_h^{k}(\lam) - \delta X_h$, until we reach the convergence with the stopping criteria we discussed before. 

From the algebraic point of view we denote with $\{E^j\}_{j=1}^{N_h}$ a base for the space $V_h$ such that we can write every element $X_h(\lam) \in V_h$ as 
\begin{equation}
\label{eq:soldecomp}
X_h(\lam) = \sum_{j=1}^{N_h} X_h^{(j)}(\lam)E^j \ , 
\end{equation}
so that we obtain the solution vector $\vec{X}_h(\lam) = \{X_h^{(j)}(\lam)\}_{j=1}^{N_h}$. We then lead back to the study of the solution $\vec{X}_h(\lam) \in \R^{N_h}$ of the system
\begin{equation}
\label{eq:lineargalerkin}
g\left(\sum_{j=1}^{N_h} X_h^{(j)}(\lam)E^j, E^i ; \lam\right) = 0 \ , \quad \forall i= 1, \dots, N_h  
\end{equation}
which, recalling the notations we introduced in Section $1.2$, corresponds to the solution of $G_h(\vec{X}_h(\lam); \lam) = 0$, where the \emph{residual vector}  $G_h$ is defined as $(G_h(\vec{X}_h(\lam); \lam))_i = g(X_h(\lam), E^i ; \lam)$.
Finally, the Newton method combined with the Galerkin finite element method, applied to our weak formulation reads: find $\delta \vec{X}_h \in \R^{N_h}$ such that 
\begin{equation}
\label{linearnewtgal}
\mathbb{J}(\vec{X}_h^k(\lam); \lam) \delta \vec{X}_h = G_h(\vec{X}_h^k(\lam); \lam) \ ,
\end{equation}
where we defined the Jacobian matrix in $\R^{N_h \times N_h}$ as 
\begin{equation}
\label{jacobiandef}
\mathbb{J}(\vec{X}_h^k(\lam); \lam))_{ij} = dg[X_h^k(\lam)](E^j, E^i; \lam) , \quad \text{for all} \quad   i,j = 1, \dots, N_h \, .
\end{equation}

%We recall that our discrete problem \eqref{eq:weakgal} is setted on a finite dimensional space $V_h$ of dimension $N_h$. 
Moreover we use standard Lagrange finite element so that $V_h = (\mathring{\mathbb{X}}_h^r)^4$ where 
\begin{equation}
 \label{spacedef}
\mathbb{X}_h^r = \{ Y_h \in C^0(\bar{\Om}) : Y_h|_K \in \mathbb{P}_r\ , \forall K \in \mathcal{T}_h\}
\end{equation}
 and 
 \begin{equation}
 \label{spacedef0}
 \mathring{\mathbb{X}}_h^r = \{ Y_h \in \mathbb{X}_h^r : Y_h|_{\p \Om} = 0\}
 \end{equation}
 is the space of globally continuous functions that are polynomials of degree $r$ on the single element of the triangulation $\mathcal{T}_h$ of the domain, which vanish on the boundary.

To provide a clear matrix representation of the application of the Galerkin method, we present the projected weak formulation. In this case the Newton method \eqref{linearnewtgal} reads: given an initial guess $X_h^0 = (u_h^0, U_h^0, \phi_h^0, \Phi_h^0) \in V_h$ for $k=0,1,\dots$ until convergence we seek $\delta X_h = (\delta u_h, \delta U_h, \delta\phi_h, \delta \Phi_h) \in V_h$ such that 
\begin{equation}
\label{eq:weakformgalerkin}
\begin{cases}
a(\delta {u_h}, w_h) + b(\delta U_h, w_h) = a(u_h^k, w_h) + b(U_h^k, w_h) \ , & \forall\,  w_h \in \mathring{\mathbb{X}}_h^r  \\
\begin{aligned}
a(\delta U_h, v_h) + c(\delta \phi_h, u_h^k, v_h) + c(\phi_h^k, \delta u_h, v_h) + \lam c(h, \delta u_h, v_h)  =  \\ a(U_h^k, v_h) + c(\phi_h^k, u_h^k, v_h) + \lam c(h, u_h^k, v_h) \ , \end{aligned} & \forall\,  v_h \in \mathring{\mathbb{X}}_h^r  \\
a(\delta {\phi_h}, \theta_h) + b(\delta \Phi_h, \theta_h) = a(\phi_h^k, \theta_h) + b(\Phi_h^k, \theta_h) \ , & \forall\,  \theta_h \in \mathring{\mathbb{X}}_h^r  \\
a(\delta \Phi_h, \psi_h) - c(\delta u_h, u_h^k, \psi_h) - c(u_h^k, \delta u_h, \psi_h) =  a(\Phi_h^k, \psi_h) - c(u_h^k, u_h^k, \psi_h) \ , & \forall\,  \psi_h \in \mathring{\mathbb{X}}_h^r  
\end{cases}
\end{equation}
and then set $X_h^{k+1} =X_h^{k} - \delta X_h$. We can finally present the matrix formulation that follows directly from \eqref{linearnewtgal} and \eqref{eq:weakformgalerkin}
\begin{equation}
\label{matrixform}
\begin{pmatrix} 
\mathbb{A}_h & \mathbb{B}_h & 0 & 0 \\
\mathbb{C}^2_h + \lam \mathbb{C}^0_h & \mathbb{A}_h & \mathbb{C}^1_h & 0 \\
0 & 0 & \mathbb{A}_h & \mathbb{B}_h \\
- \mathbb{C}^1_h - \mathbb{C}^3_h & 0 & 0 & \mathbb{A}_h 
\end{pmatrix} \begin{pmatrix} 
\delta u_h \\
\delta U_h \\
\delta \phi_h\\
\delta \Phi_h
\end{pmatrix} = \begin{pmatrix} \mathbb{A}_h u_h^k + \mathbb{B}_h U_h^k \\
\mathbb{A}_h U_h^k + \mathbb{C}^1_h u_h^k +\lam\mathbb{C}^0_h u_h^k \\
\mathbb{A}_h \phi_h^k + \mathbb{B}_h \Phi_h^k \\
\mathbb{A}_h \Phi_h^k - \mathbb{C}^1_h u_h^k \end{pmatrix} \ ,
\end{equation}
where we denoted the matrices as follows
\begin{gather*}
(\mathbb{A}_h)_{ij} = a(E^j, E^i) \ ,  \quad
(\mathbb{B}_h)_{ij} = b(E^j, E^i) \ , \quad
(\mathbb{C}^0_h)_{ij} = c(h,E^j, E^i) \ , \\
(\mathbb{C}^1_h)_{ij} = c(E^j, u_h^k, E^i) \ , \quad
(\mathbb{C}^2_h)_{ij} = c(\phi_h^k, E^j, E^i) \ , \quad
(\mathbb{C}^3_h)_{ij} = c(u_h^k, E^j, E^i) \ . 
\label{eq:matricesdef}
\end{gather*}
Note that because of the symmetry of the bracket of Monge-Amp{\'e}re, we easily obtain that it holds $\mathbb{C}^1_h \equiv \mathbb{C}^3_h$.

\subsection{Reduced Basis method}
Dealing with the approximation of a parametrized problem could be very difficult, sometimes we must therefore rely on some techniques, by which we can reduce the computational cost. With this aim in the past years many authors, to mention few works \cite{patera07:book,hesthaven2015certified,quarteroni2015reduced,morepas2017}, developed and applied the \textit{reduced order methods} (ROM), a collection of methodologies used to replace the original high dimension problem, called \textit{high fidelity approximation}, with a reduced problem that is easy to manage.

One of these methodologies is the \textit{Reduced Basis method} (RB), that consists in a projection of the high fidelity problem on a subspace of smaller dimension, constructed with some properly chosen basis functions.
%, that are exactly the solutions of the high fidelity problem once chosen a sampling technique for the parameter space. 

At the beginning this method was used for non-linear structural problems \cite{noor80:_reduc,noor83:_reduc}, but the computational complexity of this kind of equations was too big to provide a deep understanding of the bifurcation phenomena that we have analyzed.

As we have seen, the preliminary step is the projection of the weak formulation \eqref{eq:abstractweakform} in a discretized setting, which results in the Galerkin problem \eqref{eq:weakgal}. Finding a numerical solution to this problem is very challenging because of the potential high number of \textit{degrees of freedom} $N_h$.
%Instead of projecting the problem, for every new value of $\lam \in \D$, on the whole Finite Element space, 

Thus we aim at building a discrete manifold $V_N$ induced by properly chosen solution of \eqref{eq:weakgal} and then project over it. This is the description of the first step, namely the \textit{offline phase}, in which we explore the parameter space $\D$ in order to construct a basis for the reduced space of dimension $N$.
%that manifold, solving $N$ times the Galerkin problem varying $\lam$. 
On the other side, the second step, called \textit{online phase}, is the efficient and reliable part where the solutions are computed through the projection on $V_N$. This complexity reduction is based on two main key points: the assumption that it holds the affine decomposition and the fact that $N \ll N_h$. 
%The latter assumption means that we can approximate the discrete manifold with a low number of basis functions.

As in the offline phase \eqref{eq:weakgal}, given $\lam \in \D$, we seek $X_N(\lam) \in V_N$ that satisfies
\begin{equation}
\label{eq:weakred}
g(X_N(\lam),Y_N ; \lam) = 0\ , \quad \forall\ Y_N \in V_N. 
\end{equation}
At this point we have again to face with the non-linearity, and so come back to the Newton-Kantorovich  method obtaining: given an initial guess $X^0_N(\lam) \in V_N$, for every $k = 0, 1, \dots$ we seek the variation $\delta X_N \in V_h$ such that
\begin{equation}
\label{eq:weakrednewton}
dg[X_N^k(\lam)](\delta X_N, Y_N; \lam) =  g(X_N^k(\lam),Y_N ; \lam) \ , \quad \forall\ Y_N \in V_N, 
\end{equation}
and then we update the solution as $X_N^{k+1}(\lam) = X_N^{k}(\lam) - \delta X_N$ until convergence.

From the algebraic point of view and thus to do another step towards the reduced solution, we introduce the orthonormal base $\{\Sigma^m\}_{m=1}^N$ for the space $V_N$ such that we can write 
\begin{equation}
\label{eq:solreddecomp}
X_N(\lam) = \sum_{m=1}^{N} X_N^{(m)}(\lam)\Sigma^m  \ , 
\end{equation}
and denote with $\vec{X}_N(\lam) = \{X_N^{(m)}(\lam)\}_{m=1}^{N} \in \R^N$ the reduced solution vector.

Choosing properly the test element $Y_N \in V_N$ as $Y_N = \Sigma^n$ for every $1 \leq n \leq N$, we obtain the algebraic system in $\R^N$ given by
\begin{equation}
\label{eq:linearredgalerkin}
(G_N(\vec{X}_N(\lam); \lam))_n \doteq g\left(\sum_{m=1}^{N} X_N^{(m)}(\lam)\Sigma^m, \Sigma^n ; \lam\right) = 0 \ , \quad \forall n= 1, \dots, N \ ,  
\end{equation}
%We further denote with $$(G_N(\vec{X}_N(\lam); \lam))_n = g\left(\sum_{m=1}^{N} X_N^{(m)}(\lam)\Sigma^m, \Sigma^n ; \lam\right) \ ,$$ 
where we denoted with $G_N(\vec{X}_N(\lam); \lam)$ the \emph{residual reduced vector} and with $\mathbb{V}$ the transformation $N_h \times N$ matrix whose elements $(\mathbb{V})_{jm} = \Sigma^m_{(j)}$ are the nodal evaluation of the \textit{m}th basis function at the \textit{j}th node.
%the coefficients of the basis functions $\Sigma^m$ with respect to the basis functions $E^j$ of $V^h$.
Moreover, we note that \eqref{eq:linearredgalerkin} corresponds to the solution of 
\begin{equation}
\mathbb{V}^TG_h(\mathbb{V}\vec{X}_N(\lam); \lam) = 0 \ .
\label{eq:reducedproblem}
\end{equation}

Finally we can apply again the Newton method, which combined with the reduced basis method, provides the following formulation : find $\delta \vec{X}_N \in \R^{N}$ such that 
\begin{equation}
\label{linearnewtred}
\mathbb{J}_N(\vec{X}_N^k(\lam); \lam)\, \delta \vec{X}_N = G_N(\vec{X}_N^k(\lam); \lam) \ ,
\end{equation}
where $\mathbb{J}_N$ is the reduced Jacobian $\R^{N \times N}$ matrix defined as
\begin{equation}
\label{jacobianreddef}
\mathbb{J}_N(\vec{X}_N^k(\lam); \lam) = \mathbb{V}^T \mathbb{J}(\mathbb{V}\vec{X}_N^k(\lam); \lam)\mathbb{V} \  .
\end{equation}

Thus we want to construct the reduced problem through the projection on a subspace $V_N \subset V_h$ from a collection of the so called \textit{snapshots}, i.e. the solutions of the full order problem for specific values of the parameter selected by a sampling technique. The most famous strategies to construct $V_N$ are the \textit{Proper orthogonal decomposition} (POD) and the \textit{Greedy} algorithm \cite{hesthaven2015certified,patera07:book,quarteroni2015reduced}.
In this work we relied on the former, based on an ordered sampling of the interval $\D$, since its physical interpretation w.r.t the energy of the problem and because of the lack of a rigorous \textit{a posteriori error estimate} for the latter. 

Obviously POD increases the computational cost in the offline part, but simultaneously gives us a reliable representation of the reduced manifold and keep track of the energy information that we are discarding. Moreover, we can slightly modify it in order to consider multi parameter case, as the one in Section $4.3$, where the reduction shows its potentiality.
Thus, once finished the offline phase, we can build up the projection space as $V_N = span\{\Sigma^n , \ n=1,\cdots, N\}$ where $\{\Sigma^n\}_{n=1}^N$ is the basis functions set obtained through the \textit{Gram-Schmidt} orthonormalization procedure.
% starting from the set of the snapshots $\{X_N(\lam^n)\}_{n=1}^N$, and $\{\lam^n\}_{n=1}^N \subset \D$ are the values for which we solved the offline problem. In our specific computations, in order to construct properly the bifurcation plot, we used an ordered sampling of the interval $\D$. 
%Before ending this section we want to analyze in a more detailed way the system we obtained.

%In connection with the matrix representation provided previously, 
Now we show how the reduced basis method reflects the projection properties of the Galerkin method also in the online phase. Indeed the weak formulation that we obtain from the application of the Newton method at the reduced level reads as \eqref{eq:weakrednewton}, with the reduced Jacobian $\mathbb{J}_N(\vec{X}_N^k(\lam); \lam) \in \R^{N\times N}$ having the same structure of the finite element one
%: given $\lam \in \D$ and an initial guess $X_N^0 \in V_N$ for $k=0,1,\dots$ until convergence we seek $\delta X_N = (\delta u_N, \delta U_N, \delta\phi_N, \delta \Phi_N) \in V_N$ such that 
%\begin{equation*}
%\label{eq:weakformreduced}
%\begin{cases}
%a(\delta {u_N}, w_N) + b(\delta U_N, w_N) = a(u_N^k, w_N) + b(U_N^k, w_N) \ , & \forall\,  w_N \in \mathring{\mathbb{X}}_N^r  \\
%\begin{aligned}
%a(\delta U_N, v_N) + c(\delta \phi_N, u_N^k, v_N) + c(\phi_N^k, \delta u_N, v_N) + \lam c(h, \delta u_N, v_N)  =  \\ a(U_N^k, v_N)   + c(\phi_N^k, u_N^k, v_N) + \lam c(h, u_N^k, v_N)  \ ,  \end{aligned} & \forall\,  v_N \in \mathring{\mathbb{X}}_N^r \\
%a(\delta {\phi_N}, \theta_N) + b(\delta \Phi_N, \theta_N) = a(\phi_N^k, \theta_N) + b(\Phi_N^k, \theta_N) \ , & \forall\,  \theta_N \in \mathring{\mathbb{X}}_N^r  \\
%\begin{aligned}
%a(\delta \Phi_N, \psi_N) - c(\delta u_N, u_N^k, \psi_N) - c(u_N^k, \delta u_h, \psi_N) =  \\ a(\Phi_N^k, \psi_N) - c(u_N^k, u_N^k, \psi_N) \ ,  \end{aligned} & \forall\,  \psi_N \in \mathring{\mathbb{X}}_N^r  
%\end{cases}
%\end{equation*}
%and then set $X_N^{k+1} =X_N^{k} - \delta X_N$. Also in this case the matrix formulation follows directly from the considerations above and we obtain exactly the equation \eqref{linearnewtred} where
%is given by
\begin{equation}
\label{matrixformred}
\mathbb{J}_N(\vec{X}_N^k(\lam); \lam)  = \begin{pmatrix} 
\mathbb{A}_N & \mathbb{B}_N & 0 & 0 \\
\mathbb{C}^2_N + \lam \mathbb{C}^0_N & \mathbb{A}_N & \mathbb{C}^1_N & 0 \\
0 & 0 & \mathbb{A}_N & \mathbb{B}_N \\
- \mathbb{C}^1_N - \mathbb{C}^3_N & 0 & 0 & \mathbb{A}_N 
\end{pmatrix} \ ,
\end{equation}
where, if we introduce the 
%notation $\tilde{N}$ such that $N = 4\tilde{N}$ and the 
transformation matrices with respect to the different components of the solution, $\mathbb{V}_u$ and $\mathbb{V}_\phi$,  respectively for $u$ and $\phi$, we can define the reduced matrices in the following way:
\begin{gather*}
\mathbb{C}^0_N = \mathbb{V}^T\mathbb{C}^0_h\mathbb{V} \ ,  \quad
\mathbb{C}^1_N = \sum_{n=1}^{\tilde{N}}u_N^{(n)}\mathbb{V}_u^T\mathbb{C}^1_h(\Sigma^n)\mathbb{V}_u \ ,  \quad
\mathbb{C}^2_N = \sum_{n=1}^{\tilde{N}}\phi_N^{(n)}\mathbb{V}_\phi^T\mathbb{C}^2_h(\Sigma^n)\mathbb{V}_\phi \ ,  \\
\mathbb{C}^3_N = \sum_{n=1}^{\tilde{N}}\phi_N^{(n)}\mathbb{V}_\phi^T\mathbb{C}^3_h(\Sigma^n)\mathbb{V}_\phi \ , \quad 
\mathbb{A}_N = \mathbb{V}^T\mathbb{A}_h\mathbb{V} \ ,  \quad
\mathbb{B}_N = \mathbb{V}^T\mathbb{B}_h\mathbb{V}  \ .
\label{eq:matricesdefred} 
\end{gather*} 
 
Moreover, we highlight that also the reduced residual vector has the same form of the high order one, indeed it reads
 \begin{equation}
\label{residualreduced}
 G_N(\vec{X}_N^k(\lam); \lam) = \begin{pmatrix} \mathbb{A}_N u_N^k + \mathbb{B}_N U_N^k \\
\mathbb{A}_N U_N^k + \mathbb{C}^1_N u_N^k +\lam\mathbb{C}^0_N u_N^k \\
\mathbb{A}_N \phi_N^k + \mathbb{B}_N \Phi_N^k \\
\mathbb{A}_N \Phi_N^k - \mathbb{C}^1_N u_N^k \end{pmatrix} \ .
\end{equation}
 
We have illustrated the online phase, that permits an efficient evaluation of the solution and possibly related outputs for every possible choice of a different parameter $\lam \in \D$. The key point of this time saving is the so called \textit{affine decomposition}. Indeed we want the computations to be independent form the, usually very high, number of degrees of freedom $N_h$ of the \textit{true discrete problem}. In general the reduced matrices we have just presented are $\lam$-dependent and an affine-recovery technique called Empirical Interpolation Method (EIM) is needed \cite{barrault04:_empir_inter_method}. 
 
\subsection{Spectral analysis}
In the previous sections we discussed about the issue of 
%non-linearity, that lead us to the buckling model and 
the computational complexity of the problem itself, that we try to avoid using the ROM. It is clear that drawing the bifurcation diagram 
%we want to draw 
%through the continuation method
%, in order to find the bifurcation points, 
is yet a difficult task. Indeed how can we investigate the parameter space $\D$ without having any information on the position of these points? 

Taking some inspiration 
%from recent works in incompressible fluids governed by Navier-Stokes equations 
from \cite{pitton2017,pitton_quaini_2017}, where the stability property is analyzed with the help of the spectral problem,
% of the linearized problem, 
supported by the theoretical results given in Section $1.3$, we tried in this way to locate more precisely the buckling points.

%From now on we want to consider the case of the plate compressed along the edges parallel to the $y$-direction (see Figure \ref{fig:plate}), and thus assuming the compression term to be of the form $h(x,y) = -\frac{1}{2}y^2$. 
Thus we construct the eigenvalue problem for the linearized parametrized operator
%Moreover, the linearization around the null solution, provides us the following eigenproblem:

\begin{equation}
\label{eq:eigenproblem}
\begin{cases}
\Delta^2 u + \lam u_{xx} = \sigma_\lam u \ , \quad &\text{in}\ \Omega= [0,L]\times[0,1]   \\
u = \Delta u = 0 \ , \quad &\text{in}\ \p\Omega
\end{cases}
\end{equation}
where we want to find, varying the buckling parameter, the couple $(u, \sigma_\lam) \in H^1_0(\Om) \times \R$, whose components are respectively the eigenfunction and eigenvalue. %Obviously for the spectral analysis it is fundamental to specify the domain we are considering, so w
We will restrict our simulations to the square plate with $L=1$ and the rectangular one with $L=2$.
% Let us consider for a moment the eigenproblem for the square plate.

\begin{figure} [h!]
\centering   
\includegraphics[width=9cm]{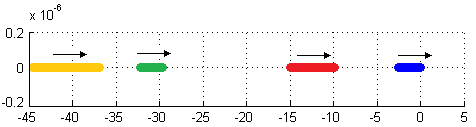}   
\caption{Behaviour of the first four eigenvalues for $\lambda \in [30,40]$.}  
\label{fig:eigenvaluefour}   
\end{figure}

We are interested in the behaviour of $\sigma_\lam$ with respect to $\lam$, in fact since the sign of the eigenvalues is strictly linked with the stability property of the solution, we aim at observing that the first eigenvalue crosses the y-axis when the plate is buckling.
This is exactly what we found, indeed in Figure \ref{fig:eigenvaluefour} for $L=1$ we can see the behaviour of the first four eigenvalues $\sigma_\lam$ for $\lam \in \left[30,40\right]$ and if we use a $\lam$-step equal to one half the crossing happens for the value $\lam = 39.5$.

This should tell us that from this point we have a change in stability properties and also the presence of a new solution branch. Indeed, as we will see, because of the symmetry, there will be at least two new branches for each simple eigenvalue. Going on with the simulations for greater values in the parameter space $\D$ we also observed the crossing of the successive eigenvalues. 

\begin{figure} [bh!]
\centering   
\includegraphics[width=9cm]{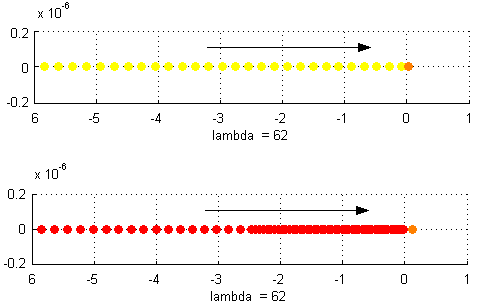}   
\caption{Double eigenvalue for the rectangular plate crossing for $\lam = 62$.}  
\label{fig:eigenvaluedouble}   
\end{figure}

Finally, if we solve the eigenvalue problem \eqref{eq:eigenproblem} for the case of $L=2$, we note that a simple eigenvalue of the square plate becomes a multiple eigenvalue with algebraic multiplicity equal to two (see Figure \ref{fig:eigenvaluedouble}). This fact has a relevant consequence from the physical point of view as we will see in the bifurcation diagram later.

What we just showed is computationally heavy to perform, so to keep in mind the efficiency as key word of the whole analysis, we tried two different ways that  validate the result and at the same time reduce the computational time.

%We asked ourself if there was a better way to individuate more precisely the buckling coefficients. Before giving an answer, let us go back to the results of Section $1.3$.
Thus we consider the linear problem \eqref{eq:linkarmsplit} but in its original form

\begin{equation}
\label{eq:eigenproblemlin}
\begin{cases}
\Delta^2 u + \lam u_{xx} = 0 \ , \quad &\text{in}\ \Omega= [0,L]\times[0,1]   \\
u = \Delta u = 0 \ , \quad &\text{in}\ \p\Omega
\end{cases}
\end{equation}
that has non trivial solutions for  $m,n = 1, 2, \dots$ given by 
\begin{equation}
\label{eq:eigenfunction}
u_{m,n} = sin\left(\frac{m\pi x}{L}\right)sin\left(n\pi y\right)  \quad 
%m,n = 1, 2, \dots
%\end{equation}
\text{if and only if} \quad
%\begin{equation}
%\label{eq:lambdavalue}
\lambda_{m,n} = \left(\frac{\pi}{L}\right)^2\left[m + \frac{n^2L^2}{m}\right]^2 , 
\end{equation}
where $u_{m,n}$ and $\lam_{m,n}$ can be considered as the eigenfunctions and eigenvalues for this new generalized eigenvalue problem. Now we have a simpler problem, indeed since $\lam$ is now the eigenvalue, there is no parametrization. This provide us also an explicit expression for the spectra, that we can use to validate the results.

Using the formula \eqref{eq:eigenfunction} we find the exact value for the eigenvalues of the problem \eqref{eq:eigenproblemlin}, that turn out to be the buckling parameter, i.e. the bifurcation point
%. In the square plate case we obtain 
\begin{equation*}
\begin{split}
L = 1:& \quad \lambda_{1,1} = 4\pi^2, \ \lambda_{2,1} = \frac{25}{4}\pi^2, \ \lambda_{3,1} = \frac{100}{9}\pi^2, \ \lambda_{4,1} = \frac{289}{16}\pi^2 ,\\
%In the very same way for the rectangular plate we find 
L = 2:& \quad \lambda_{2,1} = 4\pi^2, \ \lambda_{3,1} = \frac{169}{36}\pi^2, \ \lambda_{1,1} = \lambda_{4,1} = \frac{25}{4}\pi^2 ,
\end{split}
\end{equation*}
indeed we obtain the value $\lam_{1,1} \simeq 39.47$ predicted in Figure  \ref{fig:eigenvaluefour} for the square plate, while we note
the presence of the double eigenvalue $\lambda_{1,1} = \lambda_{4,1} \simeq 61.68$ that confirms what we saw in  Figure \ref{fig:eigenvaluedouble} for the rectangular one.

%Furthermore, we can investigate analytically the multiplicity of the eigenvalues while varying the length of the plate. We just have to impose the condition $\lambda_{m,n}=\lambda_{m+k,n}$ for some $k \in \mathbb{N}$, from which we can deduce the relation $L = \frac{\sqrt{m(m+k)}}{n}$. If we plug in the values $m=n=1$ and $k=3$ which characterize the rectangular plate, as we expected, we find the value $L=2$.

Finally, using the techniques in \cite{babuvska1991eigenvalue,MILLAR201568}, is an easy task to prove the following theorem that provides us a tool to better understand how good is our approximation.

\begin{theorem}
\label{orderofconv}
There exists a strictly positive constant $C$ such that $$|\lam - \lam_h| \leq Ch^2 ,$$ where $\lam_h$ is an approximation, dependent on the sparsity of the grid, of the true eigenvalue $\lam$.
\end{theorem}

The theorem above is crucial when we are dealing with problems for which we do not know an explicit expression of the eigenvalues. For the sake of completeness we provide in Table $1$ and Table $2$ the order of convergence results respectively for the square and rectangular plate, that agree with the theoretical ones.

\begin{table}
\centering
\begin{tabular} {cccccccc}\toprule
{L=1}  & {h = 1.e-1} & {h = 6.e-2} & {h = 1.e-2} & {h = 6.e-3} & Order & Exact \\ \toprule
$\lam_{1,1}$ & 39.91 & 39.59 & 39.48 & 39.47  & 1.98 &  39.47841\\ \midrule
$\lam_{2,1}$ & 63.70 & 62.20 & 61.70& 61.69 & 1.99 &  61.68502 \\ \midrule
$\lam_{3,1}$ & 116.63 &  111.54 & 109.73 &  109.68 & 1.97  & 109.66227 \\ \midrule
\end{tabular}
\caption{Buckling coefficients for the square plate with the average order of convergence}
\end{table}

\begin{table}
\centering
\begin{tabular} {cccccccc}\toprule
L=2	& {h = 1.e-1} & {h = 6.e-2} & {h = 1.e-2} & {h = 6.e-3} & Order & Exact \\ \toprule
$\lam_{2,1}$ & 40.74 & 39.76 & 39.48 & 39.48  & 2.05 &  39.47841\\ \midrule
$\lam_{3,1}$ & 49.15 & 46.97 & 46.35 &  46.33 & 2.34  & 46.33230 \\ \midrule
$\lam_{1,1}$ & 62.08 & 61.79 & 61.68 & 61.68 & 1.98 &  61.68502 \\ \midrule
$\lam_{4,1}$ & 67.44 & 63.02 & 61.73 & 61.69 & 2.05 &  61.68502 \\ \midrule
\end{tabular}
\caption{Buckling coefficients for the rectangular plate with the average order of convergence}
\end{table}

\begin{figure} 
\centering   
\includegraphics[width=8.1cm]{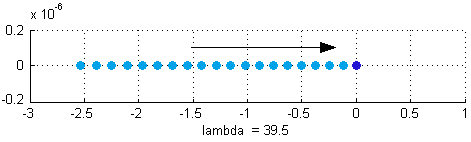}  
\caption{First eigenvalue $\lambda_{1,1}$ in the \emph{full order} case.}\label{autovall1}
\end{figure}

\begin{figure} [t!]
\centering   
\includegraphics[width=8.1cm]{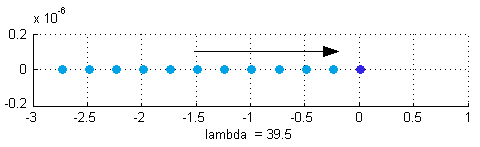}   
\caption{First eigenvalue $\lambda_{1,1}$ in the \emph{reduced order} case.}  
\label{autovall2}   
\end{figure}

To conclude this section we briefly discuss also the second straightforward way to reduce the computational complexity of solving multiple times a full order eigenvalue problem. Coming back to the parametrized eigenproblem \eqref{eq:eigenproblem}, we can apply again the Reduced Basis method, 
%in order to project the problem on a low dimensional space 
and thanks to the affine decomposition, we easily obtain in a more efficient way the same behaviour of the results discussed before, as we can see from Figures \ref{autovall1} and \ref{autovall2}. 
We do not discuss further this last approach, in fact we can embed the computation for the eigenproblem \eqref{eq:eigenproblemlin} in the offline phase.
%In fact we can compare the Figure \ref{autovall1} showing the first eigenvalue for the high order case, with the Figure \ref{autovall2} showing the same eigenvalue in the reduced order case crossing the $y$-axis.


\section{Results and test problems}
In this Section we %want to present the numerical results we obtained, that confirm all the theoretical statements discussed before. We 
will show how the buckling phenomena, i.e. the loss of uniqueness of the solution, appears in the equation through the bifurcation diagram, both in the square and rectangular plate cases.

Thus, in order to recap, let us consider the $\vK$ plate equations, with simply supported BCs, in the bi-dimensional domain $\Om = [0, L] \times [0, 1]$ given by

\begin{equation}
\label{eq:vonkarman}
\begin{cases}
\Delta^2 u + \lam u_{xx}=  \left[\phi, u\right] \ , \quad &\text{in}\ \Omega \\
\Delta^2 \phi = -\left[u, u\right] \ , \quad &\text{in}\ \Omega \\
u = \Delta u =0 \ , \quad &\text{in}\ \p\Omega \\
\phi =  \Delta \phi =0 \ , \quad &\text{in}\ \p\Omega 
\end{cases}
\end{equation}
and we are interested in the study of the solution while varying the buckling parameter $\lam$, which describes the compression along the edges parallel to the $y$-axis.

This model was previously numerically investigated by many authors \cite{refId0,chien,Reinhart1982}, but as we already said the biggest issue was the computational complexity, that we overcame by means of the Reduced Basis method.
We performed all the simulations within FEniCS \cite{LoggMardalEtAl2012a} for the full order case and RBniCS \cite{rbnics} for the reduced order one.

%In the next subsections we will consider the cases $L=1$ and $L=2$. What we expected to see, thinking about the solution of the linearized problem, is a sinusoidal behaviour of the displacement $u$. For this model it will be simple to recognize the different solutions due to the bifurcation, since they are characterized by the presence and the number of some cells in the contour plot.

%For all the simulations we chose the parameter domain as $\D = [35, 65]$, since 
The investigation done with the eigenproblem give us the necessary information that 
%of the parameter domain $\D$ With the help of the eigenproblem discussed above we know that 
the parameters responsible of the buckling live in the interval $\D = [35, 65]$, which we chose as our parameter domain.

%Finally, before discussing the application of the Reduced Basis method in order to efficiently reproduce the bifurcation diagram for the different test cases, we remark that the lack of an \textit{a posteriori error estimation} (we will explain later why) did not affect the reliability of the result we found. Indeed, we also simulated a clamped square plate under compression and forcing term $f$, properly chosen in such a way that we know the explicit expression of the solution and we can obtain the high fidelity error plot in Figure \ref{errorHF}.

%\begin{figure} [hb!]
%\centering   
%\includegraphics[width=8cm]{errorHF.png}   
%\caption{Plot of the error in $L^2(\Om)$ and $H_0^1(\Om)$ norm of the displacement $u$ for $\lam \in \D$ .}   
%\label{errorHF}   
%\end{figure}
%

\subsection{Square plate test case} 
\begin{figure} [hb!]
\centering   
\includegraphics[width=11.4cm]{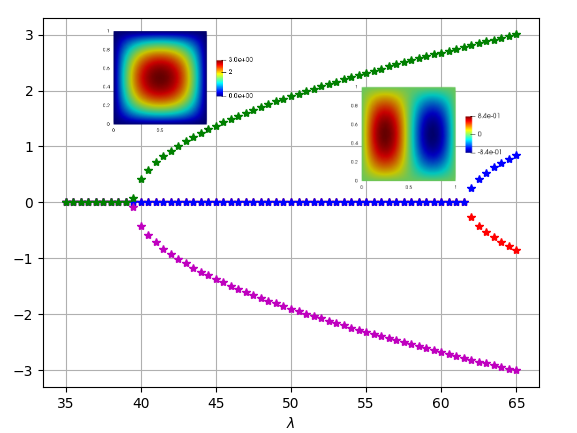}   
\caption{Bifurcation diagram for the square plate.}   
\label{bif1new}   
\end{figure}

We present the \textit{bifurcation diagram} in Figure \ref{bif1new} for the square plate case $\Om = [0,1] \times [0,1]$. The graph represents for every value of $\lam \in \D$ on the $x$-axis, the correspondent value of the full order displacement $u$ in its point of maximum modulo. How we predicted previously, we can observe the buckling phenomena from the trivial solution. Moreover, we note that the first bifurcation happens for $\lam$ value near $\lam_{1,1} \simeq 39.47$.
We did not stop at the first bifurcation, in fact choosing properly the initial guess we have been able to detect also the second bifurcation for the square plate. This result is confirming what we predicted, since for $\lam$ near $\lam_{2,1} = 61.68$ we obtain other branches.

\begin{figure}[htbp!]
        \centering
        \begin{minipage}[c]{.4\textwidth}
          \centering
\includegraphics[width=5cm]{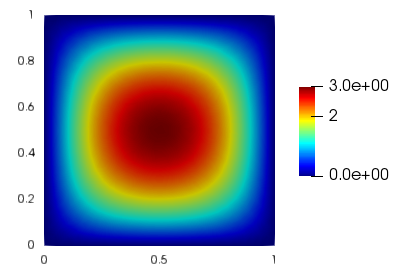}
\caption{Full order one cell solution for the displacement $u$ with $\lam = 65$ (green branch).}
\label{bifurcation1a}
        \end{minipage}%
        \hspace{20mm}%
        \begin{minipage}[c]{.4\textwidth}
          \centering
\includegraphics[width=5cm]{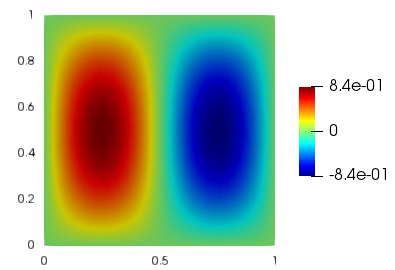}
\caption{Full order two cells solution for the displacement $u$ with $\lam = 65$ (blue branch).}
\label{bifurcation1b}
        \end{minipage}
        \end{figure}

The physical symmetry issue is evident in both buckling points and the same holds for the rectangular plate. Indeed, once we chose a bifurcation point, a solution from the upper branch is the same solution of the other one, but reflected with respect to the plate plane. Thus, for the first bifurcation near $\lam_{1,1}$, looking at the contour plot, we observe a one cell like displacement as in Figure \ref{bifurcation1a}. While if we look at the second branch, so the one near $\lam_{2,1}$, we find a two cells like displacement as in Figure \ref{bifurcation1b}.
%\begin{figure}[htbp!]
%\centering
%\includegraphics[width=8cm]{bif1_l1.png}
%\includegraphics[width=8cm]{bif2_l1.png}
%\caption{Full order one cell and two cells solutions for the displacement $u$ with $\lam = 65$ belonging to the green and blue branch, respectively on the left and on the right.}
%\label{bifurcation1}
%\end{figure}
So for the square plate case we obtained four different solutions %with respect to the same value of 
for each $\lam \geq \lam_{2,1}$.
% two of those are qualitatively different and the other ones change by the minus sign.

\begin{figure}[htbp!]
\centering
\includegraphics[width=3cm]{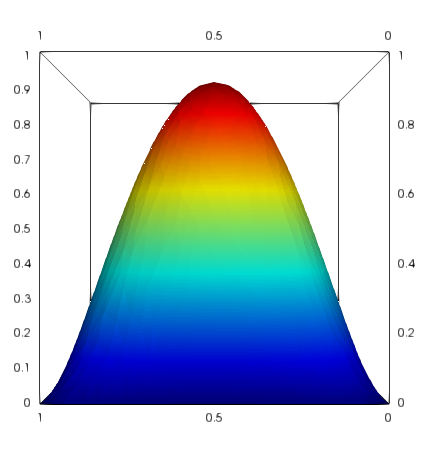}\quad\quad
\includegraphics[width=3cm]{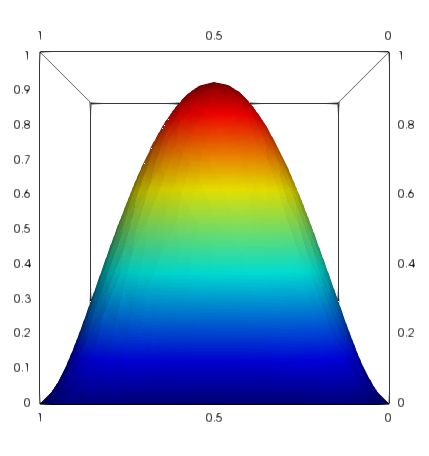} \\
\includegraphics[width=4cm]{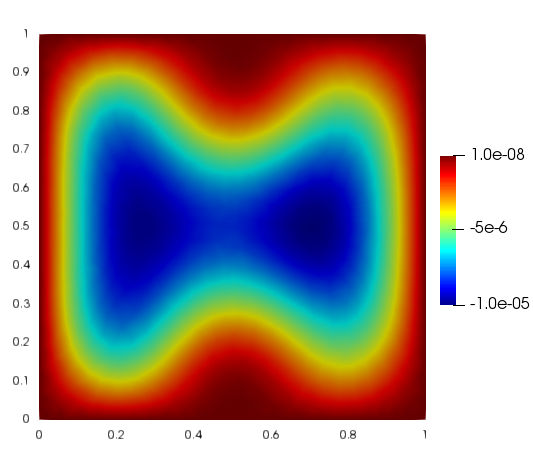}
\caption{Comparison between the full order solution (on the left) and the reduced order one (on the right) for the displacement $u$ with $\lam = 46$, belonging to the green branch. Below the reduced basis error plot.}
\label{rbsolcomp1}
\end{figure}

Moreover, the RB method worked well with this problem. Indeed as we can see in Figure \ref{rbsolcomp1}, the reduced basis solution approximates perfectly not only the behaviour but also the order of magnitude. The remarkable point is that in order to obtain the solution on the right in Figure \ref{rbsolcomp1} we just solved a linear system of dimension $5$ instead of the one given by the Galerkin full order method of order $8\cdot10^3$. 

We just saw that the Reduced Basis method provides us a useful technique, always based on a full order method (in this case is the Galerkin-Finite Element), that allows us to obtain the same results, at the cost of a small error, but with a huge amount of time saving. %The latter is mainly due to the fact that in this model the affine decomposition with respect to $\lambda$ holds. In fact in order to obtain a system that, in the online phase, is independent from the degrees of freedom of the high order problem $N_h$, we have to ensure that the equations are affine dependent by the parameter on which we are reducing the problem, and this hypothesis holds in our model.
To be more precise we present in Table $3$ a convergence results: the error between the truth approximation and the reduced one as a function of $N$. The error reported, $\mathcal{E}_N = \max_{\lam \in \D}{||u_h(\lam) - u_N(\lam)||_{H_0^1(\Om)}}$ is the maximum of the approximation error over a uniformly chosen test sample. 

We highlight that we present here just the full order bifurcation diagrams since, also in view of Table $3$ and Figure \ref{rbsolcomp1}, the reduced order one is exactly the same.

\begin{table}[h!]
\centering
\begin{tabular}{|c||c|}
    \hline
    $N$ & $\cal{E}_N$  \\
    \hline
    1 &  6.61\texttt{E}+00 \\
    2 &  6.90\texttt{E}-01 \\
 		3 &  7.81\texttt{E}-02 \\
		4 &  2.53\texttt{E}-02 \\
    5 &  1.88\texttt{E}-02 \\
		6 &  1.24\texttt{E}-02 \\
    7 &  9.02\texttt{E}-03 \\
		8 &  8.46\texttt{E}-03 \\
\hline
\end{tabular}
\label{tab:ex7_tab}
   \caption{The reduced basis convergence with respect to the number of the basis N for the square plate case.}
\end{table}
We remark that
%in literature there are many examples for the sampling methods, the most known ones are the POD and the Greedy algorithm. These methods are very useful also to give an estimate of the error $\mathcal{E}_N$. But 
we did not implemented the Greedy algorithm here, because a suitable extension of Brezzi-Rappaz-Raviart (BRR) theory for the a posteriori error estimate would be needed \cite{Brezzi1980,veroypatera2005,Grepl2007,canuto2009}. %of the solution can be obtained
%applying Brezzi-Rappaz-Raviart (BRR) theory %on the numerical approximation of non-linear problems
%\cite{Brezzi1980,veroypatera2005,Grepl2007,canuto2009}. 
However, applying BRR theory at reduced level is not straightforward and we leave it for further future investigation.

Finally, as regards computational times, a RB  evaluation $\lam \rightarrow u_N(\lam)$ requires just $t_{RB}=100$ ms for $N=8$; while the FE solution $\lam \rightarrow u_h(\lam)$ requires $t_{FE}=8.17$(s): thus our RB online evaluation is just $1.22\%$ of the FEM computational cost.

\subsection{Rectangular plate test case}
\begin{figure} [tb!]
\centering   
\includegraphics[width=12cm]{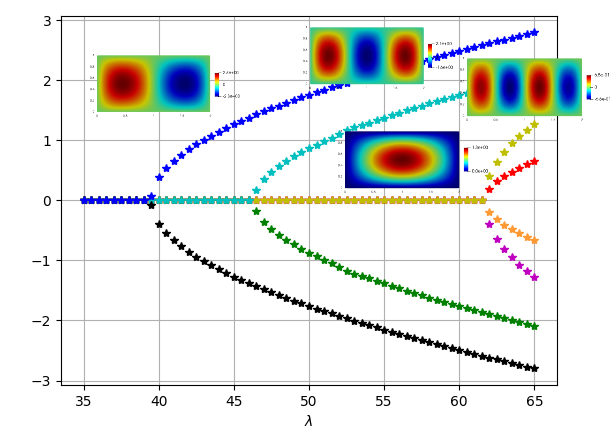}   
\caption{Bifurcation diagram for the rectangular plate.}   
\label{bif2new}   
\end{figure}

Now we analyze the case of the rectangular plate, where the domain is $\Om = [0, 2] \times [0, 1]$. A huge amount of computations led us to the bifurcation diagram in Figure \ref{bif2new}.

Now we have a different situation, in fact, varying the length of the domain, we obtained a new bifurcation and also the third one changed its properties. Note that this sensitivity with respect to the dimension of the plate is the main reason why we did not treat also $L$ as a second \textit{geometrical parameter} along with $\lam$.

\begin{figure}[htbp!]
        \centering
        \begin{minipage}[c]{.41\textwidth}
          \centering
\includegraphics[width=5cm]{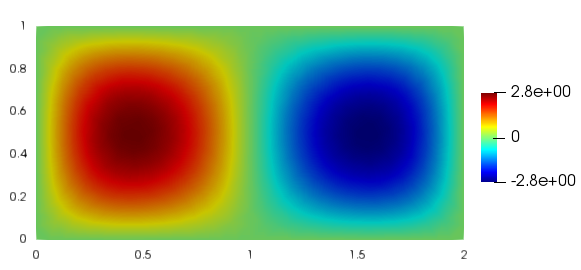}
\caption{Full order two cells solution for the displacement $u$ with $\lam=65$ (blue branch).}
\label{bifurcation2}
        \end{minipage}%
        \hspace{20mm}%
        \begin{minipage}[c]{.41\textwidth}
          \centering
\includegraphics[width=5cm]{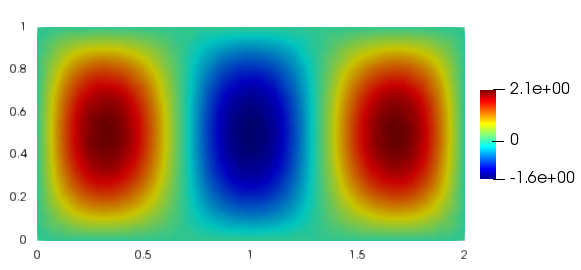}
\caption{Full order three cells solution for the displacement $u$ with $\lam=65$ (cyan branch).}
\label{bifurcation3}
        \end{minipage}
                \vspace{10mm}%
        \begin{minipage}[c]{.41\textwidth}
                  \centering
\includegraphics[width=5cm]{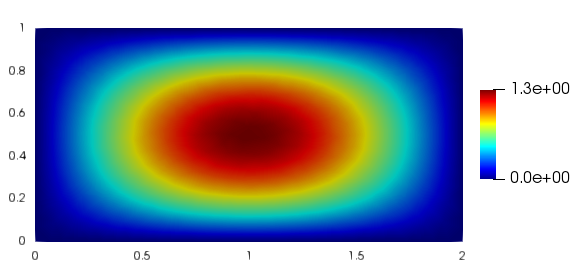}
\caption{Full order one cell solution for the displacement $u$ with $\lam=65$ (yellow branch).}
\label{bifurcation4}
        \end{minipage}
        \hspace{20mm}
                \begin{minipage}[c]{.41\textwidth}
                  \centering
\includegraphics[width=5cm]{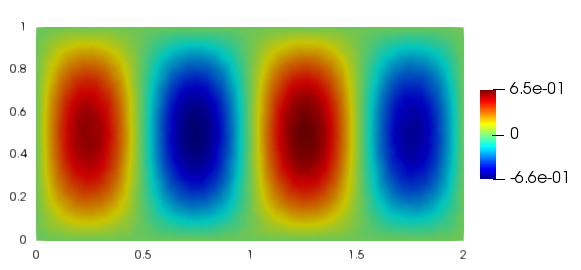}
\caption{Full order four cells solution for the displacement $u$ with $\lam=65$ (red branch).}
\label{bifurcation5}
        \end{minipage}
      \end{figure}

%We found again the parallelism between the buckling points and the eigenvalues of the generalized eigenvalue problem. 
In this case, the solution start branching, as before, from the first ordered eigenvalue $\lam_{2,1} \simeq 39.47$. %As we seen previously, we can distinguish among these solutions by their cells like displacement. 
Obviously the number of the cells in the contour plot that are formed strictly depends on the length of the domain, for example here the first bifurcation is linked with the two cell configuration as shown in Figure \ref{bifurcation2}.

 We observed also a new bifurcation for $\lam$ value near $46.5$, that is the one with three cells corresponding to the eigenvalue $\lam_{3,1} \simeq 46.33$, in Figure \ref{bifurcation3}.
 %Here the solutions are easily recognizable, as we can see in Figure \ref{bifurcation3}, by the three cells configuration.

Finally, we comment the last, qualitatively different, buckling. As we note in Figure \ref{bif2new}, as before, we have a buckling for the $\lam$ value near to $61.68$ but this time the bifurcation is linked with two eigenvalues. In fact, for the rectangular plate we have a double eigenvalue $\lam_{1,1} = \lam_{4,1}$ that is the responsible of this double bifurcation.
 In practice what we obtained is a point from which start branching two sets of different solutions with one and four cells, respectively in Figure \ref{bifurcation4} and Figure \ref{bifurcation5}. %(keeping in mind that we are always thinking about just the upper branches). 
%These solutions are again characterized by a cells like configuration that we can see respectively in Figure \ref{bifurcation4} and Figure \ref{bifurcation5} belonging to the yellow and red branches.  
%
%\begin{figure}[htbp!]
%\centering
%\includegraphics[width=5cm]{bif2_l2.png}
%\caption{Full order two cells solution for the displacement $u$ with $\lam = 65$ belonging to the blue branch.}
%\label{bifurcation2}
%\end{figure}
%
%\begin{figure}[htbp!]
%\centering
%\includegraphics[width=5cm]{bif3_l2.png}
%\caption{Full order three cells solution for the displacement $u$ with $\lam = 65$ belonging to the cyan branch.}
%\label{bifurcation3}
%\end{figure}
%
%\begin{figure}[htbp!]
%\centering
%\includegraphics[width=5cm]{bif1_l2.png}
%\caption{Full order one cell solution for the displacement $u$ with $\lam = 65$ belonging to the yellow branch.}
%\label{bifurcation4}
%\end{figure}
%
%\begin{figure}[htbp!]
%\centering
%\includegraphics[width=5cm]{bif4_l2.png}
%\caption{Full order four cells solution for the displacement $u$ with $\lam = 65$ belonging to the red branch.}
%\label{bifurcation5}
%\end{figure}

\begin{figure}[h!]
\centering
\includegraphics[width=4cm]{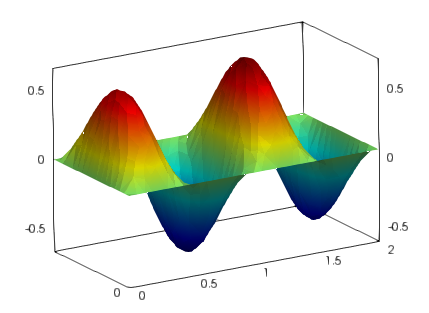}\quad\quad
\includegraphics[width=4cm]{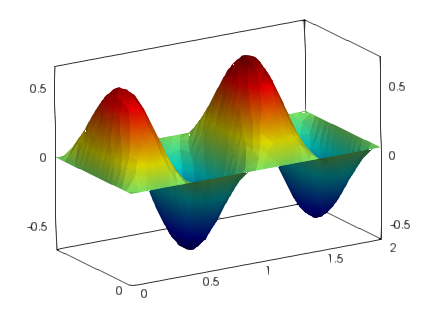}
\includegraphics[width=5cm]{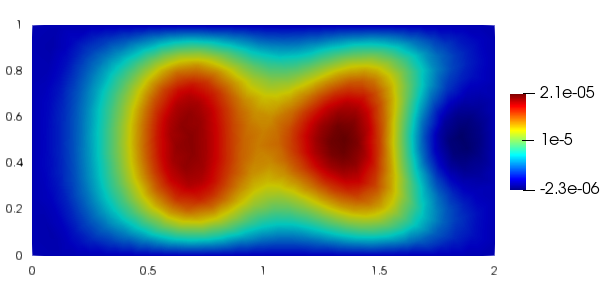}
\caption{Comparison between the full order solution (on the left) and the reduced order one (on the right) for the displacement $u$ with $\lam = 65$, belonging to the red branch. Below the reduced basis error plot.}
\label{rbsolcomp2}
\end{figure}

\vspace{-1cm}
Same conclusions regarding the convergence error $\mathcal{E}_N$ and computational savings can be established also in this case.
%Moreover, the study on the convergence error $\mathcal{E}_N$, using the Reduced Basis method, can be done also for this test case with similar results. The same holds for the computational savings that allows us to draw the bifurcation diagram efficiently (this time the computations is harder since we have to find eight branches).
Finally we show that, also for the rectangular plate, the RB method works well approximating efficiently the solution in Figure \ref{rbsolcomp2}.

\subsection{3-D bifurcation test case}
In this last numerical result section we want to extend the previous analysis in the case of two parameters, thus obtaining a 3-D bifurcation plot. Considering again the same physical phenomenon, we aim at modelling a compression along the shorter sides of the plate, which is no more uniform on that boundaries.
Interested in some practical applications of buckling plates in naval engineering, we parametrized the shape of the compression using a new parameter $\psi$.
In fact, the shape of the compression is determined by the function $h$ appearing in \eqref{eq:karm} %As an example the uniform compression over all the boundary $\p\Omega$ is modelled through the function $h = -\frac{1}{2}(x^2 + y^2)$. 
thus we can characterize the in-plane compression along $\p\Omega$ generalizing the corresponding term $\lam u_{xx}$ in \eqref{eq:eigenproblemlin} obtaining:

\begin{equation}
\label{moraeigenv}
\begin{cases}
\Delta^2 u + \lambda \,\text{div}(\sigma \nabla u) =  0 \ , \quad &\text{in}\ \Omega \\
u = \Delta u = 0 \ , \quad &\text{on}\ \p\Omega 
\end{cases}
\end{equation}  

where $\sigma: \Omega \rightarrow \R^{2\times 2}$, $\sigma \neq 0$ is the plane stress tensor field, which is assumed to satisfy the equilibrium equation:

\begin{equation}
\label{moraeqeq}
\begin{cases}
 \sigma^t = \sigma \ , \quad &\text{in}\ \Omega \\
\text{div}\,\sigma =0   \ , \quad &\text{in}\ \Omega \, .
\end{cases}
\end{equation}  

%In this way we can recover the standard compression problem, analysed in our paper, with the choice 
%$$
%\sigma = \begin{bmatrix}
%    1 & 0  \\
%    0 & 0
%\end{bmatrix} .
%$$

%Thus the compression on the whole boundary is obtained imposing $\sigma = \mathbb{I}$, where $\mathbb{I}$ is the 2-dimensional identity. 
The more general case can be studied, with the uni-axial non-uniform compression given by the stress tensor 
$$
\sigma(\psi) = \begin{bmatrix}
    \left(1- \psi\dfrac{y}{L}\right) & 0  \\
    0 & 0
\end{bmatrix} ,
$$
where $\psi \in [0, 2]$ is the parameter that takes care of the linearly varying in-plane load. Moreover, we can  recover the standard case analyzed before by choosing $\psi = 0$.

Here we want to test the strategy developed in the previous sections in this more complex case, where two parameters are involved in the bifurcation phenomenon. Here we restrict ourself to the most physically relevant behaviour, i.e. the evolution with respect to $\psi$ of the first buckling, for the generalized system 

\begin{equation}
\label{eq:vonkarmangen}
\begin{cases}
\Delta^2 u +  \lambda \,\text{div}(\sigma(\psi) \nabla u) =  \left[\phi, u\right] \ , \quad &\text{in}\ \Omega \\
\Delta^2 \phi = -\left[u, u\right] \ , \quad &\text{in}\ \Omega \\
u = \Delta u =0 \ , \quad &\text{in}\ \p\Omega \\
\phi =  \Delta \phi =0 \ , \quad &\text{in}\ \p\Omega\  .
\end{cases}
\end{equation}

\begin{figure}[htbp!]
%        \centering
        \begin{minipage}[c]{.45\textwidth}
%          \centering
%\hspace{-10mm}%
\includegraphics[width=6cm]{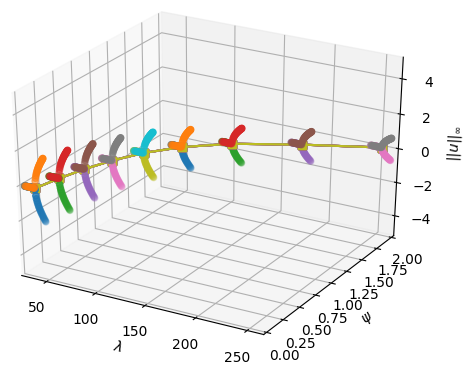}   
\caption{3-D bifurcation diagram for the square plate.}   
\label{fig:3D_bifurcation_plot}   
        \end{minipage}%
        %\hspace{10mm}%
        \hfill
        \begin{minipage}[c]{.45\textwidth}
%          \centering
\includegraphics[width=5.2cm]{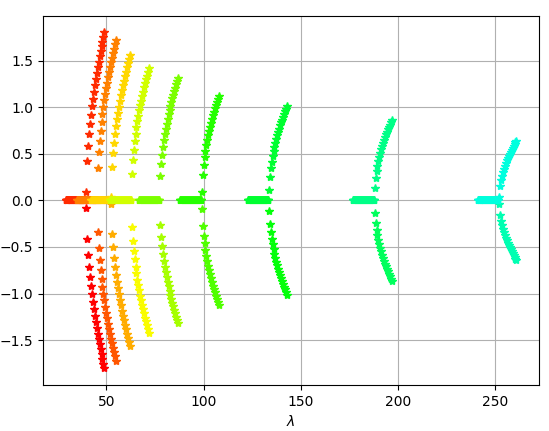}   
\caption{2-D projected bifurcation diagram for the square plate.}   
\label{fig:HO_bifurcartion_diagram}   
        \end{minipage}
        \end{figure}

Now we can show some preliminary results on the behaviour of the first buckling for the system \eqref{eq:vonkarmangen}.
First of all we can observe in the Figure \ref{fig:3D_bifurcation_plot} the 3-D bifurcation plot for the square plate, in which we are describing the first bifurcation point and the post-buckling behaviour, without loss of generality,  for each uniformly sampled $\psi \in [0, 2]$. For the sake of clearness we show also in Figure \ref{fig:HO_bifurcartion_diagram} the 2-D version of the 3-D plot just presented.

%\begin{figure} [htb!]
%\centering   
%\includegraphics[width=6cm]{3D_bifurcation_plot.png}   
%\caption{3-D bifurcation diagram for the square plate.}   
%\label{fig:3D_bifurcation_plot}   
%\end{figure} 
%
%\begin{figure} [htb!]
%\centering   
%\includegraphics[width=6cm]{HO_bifurcartion_diagram.png}   
%\caption{2-D projected bifurcation diagram for the square plate.}   
%\label{fig:HO_bifurcartion_diagram}   
%\end{figure} 

As we can see, the methodology presented in the previous section was able to well detect in the reduced phase the first bifurcation points with respect to the new parameter $\psi$ just introduced. Moreover, we were able to capture correctly the post buckling behaviour, with results validated by the former analysis with $\psi = 0$.

Finally we want to remark that here the necessity of the Reduced Order Models is still more evident. In fact, considering only the full order problem, we have to solve a huge linear system (as many times as the following nested iteration):  for each Newton step, for each $\lambda$ in the parameter domain, for each (selected) $\psi$ and for each initial guess if one is interested on multiple branches.

\section{Future perspectives and developments}
In this work we have presented a methodology to properly detect a bifurcation phenomena  at different levels, with the support of strong and consolidated theoretical results and the help of computational reduction strategies, allowing to predict efficiently the buckling. We showed the connection of these physical phenomena with the eigenvalue problem, which is fundamental in dealing with different geometry or more complex applications. 
Several numerical tests were performed confirming the strength of the reduced basis method and its reliability, also in non-linear context.
The recovery of eight of the possible solutions shows that, approaching with complex non-linear problems, we need to rely on some backup tool in order to verify that the solution we found is the one we are interested in.
Moreover, with this work we showed the consistency of the theoretical results with the numerical ones, but also the necessity to investigate the reduction strategies to apply this methodology on more complex, real applications. The extension of the results for the multi parameter application could be also more relevant since the increasing computational cost.

We plan to extend this work in different directions. The first one is towards the Brezzi-Rappaz-Raviart theory providing the model with an ``a posteriori error estimate". Furthermore, we want to apply this methodology to other kind of problems, such as in fluid structure interaction models, as well as in vibro-acoustics and fluid mechanics frameworks. From the continuum mechanics point of view, we are also interested in the study of different type of plate models, such as the Saint Venant-Kirchhoff, as well as the extension toward the three dimensional $\vK$ equations.

\section*{\ackname}
 The authors thank Dr.~F.~Ballarin (SISSA) for his great help with the RBniCS software and precious discussion. The authors thank Prof.~A.~T.~Patera for the inspiring conversations and valuable time. This work was supported by European Union Funding for Research and Innovation through the European Research Council (project H2020 ERC CoG 2015 AROMA-CFD grant 681447, P.I. Prof. G. Rozza) by the INDAM-GNCS 2017 project ``Advanced numerical methods combined with computational reduction techniques for parameterised PDEs and applications", and by the MIT-FVG project ROM2S "Reduced Order Methods at MIT and SISSA".
%%-----------------------------
\bibliography{reference}
\bibliographystyle{plain}
%%-----------------------------

%\begin{acknowledgements}
%If you'd like to thank anyone, place your comments here
%and remove the percent signs.
%\end{acknowledgements}

% BibTeX users please use one of
%\bibliographystyle{spbasic}      % basic style, author-year citations
%\bibliographystyle{spmpsci}      % mathematics and physical sciences
%\bibliographystyle{spphys}       % APS-like style for physics
%\bibliography{}   % name your BibTeX data base

%% Non-BibTeX users please use
%\begin{thebibliography}{}
%%
%% and use \bibitem to create references. Consult the Instructions
%% for authors for reference list style.
%%
%\bibitem{RefJ}
%% Format for Journal Reference
%Author, Article title, Journal, Volume, page numbers (year)
%% Format for books
%\bibitem{RefB}
%Author, Book title, page numbers. Publisher, place (year)
%% etc
%\end{thebibliography}

\end{document}